\documentclass[a4paper,12pt]{article}

\usepackage{amsthm}
\usepackage{amsfonts}
\usepackage{amsmath}
\usepackage{times}

\newcommand{\prob}[1]{\mathsf{\Pr}\left(#1\right)}
\newcommand{\EXP}[1]{\mathsf{E}\left(#1\right)}
\newcommand{\VAR}[1]{\mathsf{Var}\left(#1\right)}
\newcommand{\nth}{^{\textrm{th}}}
\newcommand{\remove}[1]{}
\newtheorem{theorem}{Theorem}
\newtheorem{lemma}{Lemma}
\newtheorem{definition}{Definition}
\newtheorem{remark}{Remark}
\newtheorem{property}{Property}

\newcommand{\NH}{N\hspace{-3pt}H}

\def\rar{\rightarrow}

\def\ep{\epsilon}
\def\lam{\lambda}
\def\al{\alpha}
\def\del{\delta}

\def\Th{\Theta}

\begin{document}

\title{On the Topological Properties of the One Dimensional
  Exponential Random Geometric Graph}

\author{
  Bhupendra Gupta\thanks{Department of Mathematics and Statistics, Indian
    Institute of Technology, Kanpur 208016 INDIA. email:
    \texttt{bhupen@iitk.ac.in}.}, Srikanth K.~Iyer\thanks{Department of
    Mathematics, Indian Institute of Science, Bangalore 560012
    INDIA. email: \texttt{skiyer@math.iisc.ernet.in}} and
  D.~Manjunath\thanks{Department of Electrical Engineering, Indian
    Institute of Technology, Bombay, Mumbai 400076,
    INDIA. email:\texttt{dmanju@ee.iitb.ac.in}.
  }
}

\date{}

\maketitle

\begin{abstract}
  In this paper we study the one dimensional random geometric (random
  interval) graph when the location of the nodes are independent and
  exponentially distributed. We derive exact results and limit
  theorems for the connectivity and other properties associated with
  this random graph. We show that the asymptotic properties of a graph
  with a truncated exponential distribution can be obtained using the
  exponential random geometric graph.
\end{abstract}


\newpage

\section{Introduction} \label{sec:intro}

We consider random geometric graphs (RGGs) in one dimension,
$G_n(\lambda,r)$, with vertex set $V_n=\{X_1,\ldots, X_n\}$ and
edge set $E=\{(X_i,X_j): |X_i-X_j| \leq r\}$, where $X_i$ are
i.i.d. exponential with mean $\lambda^{-1}$ and $r$ is called the
``cutoff range''. Here $X_i$ is used to denote the $i\nth$ vertex
and its location. $G_n(\lambda,r)$ will be called an exponential
RGG. We derive formulas and recursive algorithms when the number
of nodes $n$ and cutoff $r$ are fixed. We then derive asymptotic
results for the probability of connectivity, and weak law results
for the number of components, total uncovered area etc. Strong
law asymptotics are derived for the connectivity and largest
nearest neighbor distances. We also obtain strong law results
when the $X_i$ are i.i.d. truncated exponential.

\subsection{Previous Work and Background}
The topological properties of RGGs have applications in wireless
communication and sensor networks (e.g., \cite{Gupta98}), cluster
analysis (e.g., \cite{Godehardt90,Godehardt93}), classification
problems in archaeological findings, traffic light phasing, and
geological problems (e.g., \cite{Godehardt96}), and also in their
own right (e.g., \cite{Penrose03}).

The following are some results motivated by random wireless
networks. For a network of $n$ nodes distributed unformly inside
the unit circle, \cite{Gupta98} obtains the asymptotic threshold
function for the critical transmission range to be
$\sqrt{\frac{\log n}{n}}$. More precisely, they show that with
cutoff $r_n = \sqrt{\frac{\log n + c_n}{n}},$ the graph on $n$
uniform points in the unit circle will be connected with
probability approaching one iff $c_n \rar \infty$. A result that
enables the nodes to control local topological properties and
work towards obtaining a connected network is derived in
\cite{Xue02}.  Here it is shown that for a static network with
$n$ nodes uniformly distributed over the unit circle, if each
node is connected to $( 5.1774 \log n) $ nodes, then the network
is asymptotically connected.  This problem has also been studied
empirically in the context of multi-hop slotted Aloha networks
\cite{Kleinrock78}. The transmission radius for connectivity of a
network when the placement of the nodes follows a Poisson process
in dimensions $d \leq 2$ is derived in \cite{Cheng89}.

The following is a sample of the results from the study of RGGs in
their own right. When $n$ nodes are uniformly distributed in the
$d$-dimensional unit cube, the following is shown in \cite{Penrose99}
for any $l_p$ metric.  Start with isolated points and keep adding
edges in order of increasing length.  Then, with a very high
probability, the resulting graph becomes $(k+1)$-connected at the same
edge length $r^*(n)$ at which the minimum degree of the graph becomes
$k+1$, for $k \geq 0$.  With $k=0$, this result means that the graph
becomes connected with high probability at the same time that the
isolated vertices disappear from the graph.  \cite{Appel97a,Appel97b}
is a similar study for the $l_\infty$ norm. The best introduction to
the study of RGGs via their asymptotic properties is \cite{Penrose03}.

Observe that the results cited above are all for the asymptotic case
with $n \to \infty$. Exact analysis of finite networks is important
because the asymptotes may be approached very slowly. Exact analysis
of finite networks have been considered in
\cite{Desai02,Desai05,Godehardt96}.  The probability of a connected
network when $n$ nodes are uniformly distributed in $[0,1]$, with all
nodes having the same transmission range $r$ was derived in
\cite{Desai02,Godehardt96}. In \cite{Godehardt96}, algorithms for
various connectivity properties of a one dimensional RGG with uniform
distribution of nodes on the unit interval are derived.  The
probability of a specified labeled subgraph with edge set $E=\{
(X_i,X_j): |X_i - X_j| \leq w_{i,j} \}$ was obtained in
\cite{Desai05}. This is then used to calculate exact probabilities for
many network properties. This generalization dispenses with the
requirement that the cutoff range be the same for all node pairs.

While the results described above are all for the case of an RGG when
the nodes are distributed uniformly in a finite operational area, they
can be extended to RGGs where the density of the node locations is
arbitrary but has bounded support. The asymptotic behaviour here is
similar to that of a graph with uniform distribution of nodes
\cite{Penrose03}.

If the region of deployment is large, it makes sense to consider
distributions with unbounded supports. As in other applications, this
would offer us a wide variety of nice distributions which can be used
to answer many interesting questions regarding the RGGs. Most
interesting results for such densities depend on the tail behavior of
the underlying distribution.

In this paper we primarily consider RGGs where the distribution
of the node locations are i.i.d. exponential. The motivation is
from random wireless sensor networks. Consider the deployment of
intrusion detection sensors along a border. The cost of the
sensors is expected to be significantly less than the cost of a
`regular' deployment.  We remark here that a class of such
relatively inexpensive devices called \emph{smart dust}
\cite{url-smartdust} are actually available! Hence, it is not
unreasonable to expect that the sensors will be deployed by a
random dispersion onto the border line.  If the point from where
they are dispersed is treated as the origin, it is reasonable to
expect that the distribution of the sensor nodes will be dense
near the origin and sparse away from it. Thus it is important to
consider non-uniform distribution of the nodes. Further, analysis
of networks with a finite number of nodes would also be very
useful.

We remark here that the asymptotic results that are governed by the
clustering of the nodes near the mode, e.g., maximum vertex degree,
are obtained as in the case of RGGs with finite support
\cite{Penrose03}.  In contrast, characteristics such as the largest
nearest-neighbor distance, connectivity distance, minimum vertex
degree etc. for densities with unbounded support are dependent on the
tail-behavior and connectivity distances for normally distributed
nodes are obtained in \cite{Penrose03}.

\subsection{Summary of Results and Outline of Paper}
Consider the exponential RGG, $G_n(\lambda,r)$ with node locations
$\{X_1, \ldots, X_i, \ldots, X_n\}.$ Let $X_{(i)}$, denote the
distance of the $i\nth$ node from the origin or the $i\nth$ order
statistics of the random sample $\{X_i\}$. Let $X_{0} = 0$ and define
$Y_{i} := X_{(i+1)} - X_{(i)}$ $i=0,1, \ldots (n-1).$ The following is
a key result that we will use quite often in the remainder of this
paper. From \cite{Arnold93}, we have the following lemma.
\begin{lemma}
  $Y_1,Y_2, \ldots ,Y_{n-2}, Y_{n-1}$ are independent exponential
  random and the means are $((n-1)\lambda)^{-1}, ((n-2)\lambda)^{-1},
  \ldots (2\lambda)^{-1}, (\lambda)^{-1}$ respectively.
  \label{lemma:independent-yi}
\end{lemma}

The lemma follows from the fact that the minimum of $m$ i.i.d.
exponentials of mean $1/\lambda$ is an exponential of mean $(m
\lambda)^{-1}$ and from the memoryless property of the exponential
distribution.

The rest of the paper is organized as follows.  In
Section~\ref{sec:connectivity}, we derive the exact expression
for the probability of connectivity $P^c_n$ of the one
dimensional exponential RGG with $n$ nodes. In
Theorem~\ref{thm:conn-asymptote} we show that $P_n^c \rar P_c$ as
$n \rar \infty$, where $0 < P_c < 1.$ This limit and all other
asymptotics hold under the condition that $\lambda r$ is fixed or
converges to a constant.  This is in contrast to the limiting
results for the uniform and the normal case where the limiting
results under the condition that $r_n \to 0$ (see
\cite{Penrose03}). In Section~\ref{sec:components}, we first give
a recursive formula for the distribution of the number of
components for finite $n$. In Theorem~\ref{thm:comp-asymptote} we
show that this distribution converges as $n \rar \infty$ and in
Theorem~\ref{thm:dist_size_m_comp} we obtain limiting
distribution for the number of components of size $m.$
Section~\ref{sec:redundant} provides a recursive formula for
computing the distribution of the number of redundant nodes, nodes
that can be removed without changing the connectivity of the
network. In Section~\ref{sec:expected-node-degree} we
characterize the degree of a node by obtaining the asymptotic
expectation of the degree in
Theorem~\ref{thm:expected-degree-asym}.
Section~\ref{sec:span-asymptote} deals with the span and the
uncovered part of the network. In
Theorem~\ref{thm:holes-convergence}, we show that the span of the
network converges to $\infty$ with probability $1.$ However, the
total number of holes (gaps between ordered nodes of length
greater than $r$) and the total length of the holes converge in
distribution. An interesting upshot of this result is that though
the span of the network diverges, the probability of connectivity
converges to a non-zero constant. Thus we can achieve (by taking
$n$ large) an arbitrarily large coverage with high probability,
without diminishing the probability of connectivity.
Theorem~\ref{thm:span-distribution-asym} derives the asymptotic
distribution of the span of the network.

In Section~\ref{sec:strong-laws}, we derive strong law results for
connectivity and largest nearest neighbor distances in
Theorem~\ref{thm:slln_exp}. Finally, in
Section~\ref{sec:truncated-exponential} we consider RGGs where the
node locations are drawn from a truncated exponential distribution,
i.e., the exponential restricted to $(0,T)$. show that the asymptotic
results for the truncated exponential RGG can be derived using
properties of the exponential RGG $G_n(\lambda,r)$. We first define
monotone properties and the strong and weak thresholds for the cutoff
distance $r$ for monotone properties. In
Theorem~\ref{thm:threshold_equivalence_truncated_exp} we show the
equivalence of strong and weak thresholds for monotone properties in a
truncated exponential RGG and an RGG constructed by considering the
first $n$ nodes of an exponential RGG. Using this, in
Theorem~\ref{thm:threshold_connectivity_truncated_exp} we obtain the
cutoff thresholds for the RGG to be connected.
Theorem~\ref{thm:slln_trunc_exp} obtains the strong law for the
connectivity and largest nearest neighbor distances.

We remark here that many of the results that we derive for the one
dimensional exponential network can also be extended to the case of
the nodes being distributed according to the double exponential
distribution which is just the exponential density defined on the
entire real line. It has the density $\frac{\lambda}{2}e^{-\lambda
  |x|}$ for $-\infty < x < \infty$. We will derive only the
probability of connectivity for the double exponential case.

\section{Connectivity Properties} \label{sec:connectivity} Let $P^c_n$
denote the probability that a network of $n$ nodes each with a
transmission range $r$ is connected.  For the network to be connected
we must have $Y_i=X_{(i+1)} - X_{(i)} \leq r,$ $\forall \; i =
1,2,\ldots,(n-1)$. From Lemma~\ref{lemma:independent-yi}, the
following is straightforward.
\begin{property}
  $P^c_n$ is given by
  \begin{eqnarray}
    P^c_n \ \ = \ \ \prod_{i=1}^{n-1} \prob{Y_i \leq r} \ \ = \ \
    \prod_{i=1}^{n-1} (1-e^{-(n-i){\lambda}r}) = \prod_{i=1}^{n-1}
    (1-e^{-i \lambda r}).
    \label{eq:prob_connect_exp}
  \end{eqnarray}
  \label{prop:P^c_n}
\end{property}

We now derive the probability that a network constructed using the
double exponential distribution is connected.  We condition on the
event that of the $n$ nodes, $k$ nodes are in $(0,\infty)$ and $n-k$
are in $(-\infty,0)$. Label the positive observations as $U_i,$
$i=1,\ldots,k,$ and the absolute values of the negative observations
as $V_i,$ $i=1,\ldots,(n-k)$.  Then the $U_i$ and $V_i$ are
independent exponential variables with mean $1/\lambda$.  If the
network of $U$ values is connected and the network of $V$ values is
connected and the distance between the $ U_{(1)}$ and $-V_{(1)}$ is
less than $r$, then the network will be connected. Note that from
Lemma~\ref{lemma:independent-yi}, it follows that $U_{(1)}$ and
$V_{(1)}$ are independent of whether the networks on the positive and
negative halves are connected or not. Thus, the probability that the
network is connected, $P^c_n(D)$, will be
\begin{equation}
  \!\! P^c_n(D) \!\! = \!\! \sum_{k=1}^{n-1} \!\!
  \binom{n}{k} (1/2)^n  \prob{U_{(1)}+
    V_{(1)} \leq r \mid X_{(k)} < 0, X_{(k+1)} > 0 }  P^c_k  P^c_{n-k}
  + \frac{P^c_n}{2^{n-1}}.
  \label{eq:prob_connect_double_exp}
\end{equation}
The densities of $U_{(1)}$ and $V_{(1)}$ conditioned on the event
$\{ X_{(k)} < 0, X_{(k+1)} > 0 \}$ will be
\begin{eqnarray*}
  f_{U_{(1)}}(u) & = & k \lambda e^{-k\lambda u} \qquad 0<u<\infty, \\
  f_{V_{(1)}}(v) & = & (n-k) \lambda e^{-(n-k)\lambda v} \qquad
  - \infty < v < 0.
\end{eqnarray*}
The density of $\left(U_{(1)}+V_{(1)}\right)$,
$g_{U_{(1)}+V_{(1)}}(z)$, and hence the probability that
$U_{(1)}$ and $V_{(1)}$ are connected, is now straightforward;
\[ g_{ U_{(1)}+V_{(1)}}(z) = \begin{cases}
  \frac{k(n-k)\lambda}{n-2k} \left(e^{-k\lambda z} -
    e^{-(n-k)\lambda z} \right)  &
  \hspace{0.1in} \mbox{ if $2k \neq n$}\\
  (k \lambda)^2 z e^{-k \lambda z} &
  \hspace{0.1in} \mbox{ if $2k=n$},
\end{cases}
\]
and
\begin{equation}
  \prob{U_{(1)}+ V_{(1)}\leq r} = \begin{cases} 1 +
    \frac{1}{n-2k}\left( k e^{-(n-k)\lambda r} - (n-k)e^{-k \lambda r}
    \right) & \hspace{0.1in} \mbox{if $2k \neq n$}\\ 1-e^{-k \lambda
    r}(1+k \lambda r) & \hspace{0.1in} \mbox{if $2k=n$}.
  \end{cases}
  \label{eq:prob_UV_connection}
\end{equation}
%

%
Using (\ref{prop:P^c_n}) and (\ref{eq:prob_UV_connection}) in
(\ref{eq:prob_connect_double_exp}) we obtain the following.
\begin{property}
  If the $X_i$ are i.i.d. double exponential with zero mean, then the
  probability that the network is connected, $P^c_n(D)$, is given by
  \begin{eqnarray}
    P^c_n(D) & = & \frac{1}{2^n} \sum_{{k=0} \atop{k \neq
        n/2}}^{n} \binom{n}{k} \ P^c_k \ P^c_{n-k} \left( 1 +
      \frac{1}{n-2k}\left(k  e^{-(n-k)\lambda r} - (n-k) e^{-k\lambda
          r} \right) \right)  \nonumber \\
    && \hspace{-0.2in} \ + \ \frac{(P^c_{n/2})^2}{2^n} (1-e^{-n \lambda
        r/2}(1+n \lambda r/2)).
    \label{eq:prob_connect_double_exp_formula}
  \end{eqnarray}
  \label{prop:P^c_n(D)}
\end{property}
In (\ref{eq:prob_connect_double_exp_formula}), we have defined
$P^c_0=1$. Also, the last term will be necessary only when $n$ is
even.
\begin{theorem}
Let $P^c_n$ and $P^c_n(D)$ denote the probability that the
exponential and double exponential random geometric graphs
respectively, with $n$ vertices, parameter $\lam,$ and cutoff $r$
are connected. Then, for some real number $P_c$, $0 < P^c < 1,$
  \begin{enumerate}
  \item $\lim_{n \to \infty} P^c_n = P^c ,$
  \item $\lim_{n \to \infty} P^c_n(D) = (P^c)^2$
  \end{enumerate}
  \label{thm:conn-asymptote}
\end{theorem}
\proof{

  Consider the first part of the theorem. Taking logarithms on both
  sides of (\ref{eq:prob_connect_exp}) we get
  \begin{eqnarray}
    \lim_{n \to \infty} \ln(P^c_n) & = & \sum_{i=1}^{\infty}
    \ln(1-e^{-i \lambda r}) = \sum_{i=1}^{\infty}
    \sum_{j=1}^\infty \ \frac{-(e^{-i \lambda r})^j}{j} \nonumber \\
    & = & - \sum_{j=1}^{\infty} \frac{1}{j} \sum_{i=1}^\infty \ (e^{- j
      \lambda r})^i = - \sum_{j=1}^{\infty} \frac{1}{j}  \frac{e^{- j \lambda
        r}}{1-e^{- j \lambda r}}.
    \label{eq:sum-for-C}
  \end{eqnarray}
  Applying the ratio test we see that the series converges to a finite
  value $\ln P^c$. Since $-\infty < \ln(P^c):=\lim_{n \to \infty} \ln(P^c_n) < 0$
  we get $0 < P^c < 1$.

  Now consider the second part of the theorem statement. Let $L_n$ be
  the number of nodes to the left of the origin when $n$ nodes are
  distributed on the real line. By the strong law of large numbers,
  $\frac{L_n}{n} \stackrel{\mathrm{a.s}}{\rightarrow} \frac{1}{2}$.
  This implies that for any $\epsilon > 0$, there exists a finite
  $m(\epsilon)$ such that
  \begin{equation}
    \prob{\sup_{n \geq m(\epsilon)}
      \left|L_n - \frac{n}{2}\right| > n \epsilon} <  \epsilon.
    \label{eq:slln-m-eps}
  \end{equation}
  To make the notation below simpler, we will assume that $n$ is odd.
  Let $n > m(\epsilon)$.  Define
  \begin{eqnarray}
    A_{n,k} \!\! &:=& \! \! \left(1 + \frac{1}{n-2k} \left(k
        e^{-(n-k)\lambda r} - (n-k) e^{-k\lambda r} \right) \right),
    \; k=1,\ldots,(n-1). \nonumber \\
  \end{eqnarray}
  Using the preceding definition for $A_{n,k}$, we can write
  (\ref{eq:prob_connect_double_exp_formula}) as
  \begin{eqnarray*}
    P^c_n(D) &=& \sum_{k=1}^{n-1} \binom{n}{k} \frac{1}{2^n} P^c_k
    P^c_{n-k} A_{n,k} + \frac{P^c_n}{2^{n-1}} \\
    &=& \sum_{k: \ |k-n/2|  \leq n \epsilon}^{n-1} \binom{n}{k}
    \ \frac{1}{2^n} P^c_k P^c_{n-k} A_{n,k} +
    \sum_{k:\ |k-n/2| > n \epsilon }^{n-1}
    \binom{n}{k} \frac{1}{2^n} P^c_k P^c_{n-k} A_{n,k} +
    \frac{P^c_n}{2^{n-1}}.
  \end{eqnarray*}
  We first consider the second term in the above equation;
  \begin{equation}
    \sum_{|k-n/2| > n \epsilon
      }^{n-1} \!\! \binom{n}{k} \frac{1}{2^n} P^c_k \ P^c_{n-k} \ A_{n,k}
    \leq \sum_{|k-n/2| > n \epsilon}^{n-1} \!\!\binom{n}{k}
    \frac{1}{2^n} <  \epsilon.
    \label{eq:2side-ineq0}
  \end{equation}
  The last inequality is derived by first observing that since $k$ is
  the number of nodes to the left of the origin, the summation
  corresponds to the probability of $\{|L_n -n/2| > n \epsilon\}$ and
  then applying (\ref{eq:slln-m-eps}). Now consider the first sum,
  \begin{equation}
    \sum_{|k-n/2| \leq n \epsilon}^{n-1} \binom{n}{k}
    \frac{1}{2^n} P^c_k P^c_{n-k} A_{n,k} \leq
    \left(P^c_{n/2-n\epsilon}\right)^2 \sum_{{k=1}\atop{|k-n/2| \leq n
        \epsilon}}^{n-1} \binom{n}{k} \frac{1}{2^n} <
    \left(P^c_{n/2-n\epsilon}\right)^2.
    \label{eq:2side-ineq1}
  \end{equation}
  The first inequality is true because $P^c_{n}$ is decreasing in $n$
  and $A_{n,k} \leq 1$ ($A_{n,k}$ is a probability). The last
  inequality is true since the sum is less than 1. Also, note that
  $\lim_{n \to \infty} P^c_{n/2-n\epsilon} = P_c$.

  We can also write the following inequality.
  \begin{equation}
    \! \sum_{|k-n/2| \leq n \epsilon}^{n-1} \binom{n}{k}
    \frac{1}{2^n} P^c_k \ P^c_{n-k} \ A_{n,k} \geq
    \left(P^c_{n/2+n\epsilon}\right)^2 \
    \sum_{|k-n/2| \leq  n \epsilon}^{n-1} \binom{n}{k}
    \frac{1}{2^n} A_{n,k}  \! \geq \! P_c^2 (1-\epsilon)^2.
    \label{eq:2side-ineq2}
  \end{equation}
  The first inequality is true because $P^c_n$ is decreasing in $n$.
  To see why the second inequality is true, we first note that
  $\lim_{n,k\to \infty} A_{n,k} =1$. Hence for large $|k-n/2| > n
  \epsilon $ and large $n$ $A_{n,k} > (1-\epsilon)$. Combining this
  observation with (\ref{eq:slln-m-eps}) and noting that $P^c_n$
  converges monotonically to $P_c$ we can write the second inequality in
  (\ref{eq:2side-ineq2}). Thus, from (\ref{eq:2side-ineq1}) and
  (\ref{eq:2side-ineq2}), we get
  \begin{displaymath}
    \lim_{n \to \infty} \ \sum_{|k-n/2| \leq n
      \epsilon}^{n-1} \binom{n}{k} \  \frac{1}{2^n} P^c_k \ P^c_{n-k} \
    A_{n,k}  \ = \ \left(P_c\right)^2.
  \end{displaymath}
  Combining this with (\ref{eq:2side-ineq0}), the second part of the
  theorem is proved.  \hfill $\Box$

  }

Numerical evaluation shows that both $P^c_n$ and $P^c_n(D)$ converge
rapidly.

\section{Components in the Network}
\label{sec:components}

A sequence of connected nodes which are followed and preceded by a
disconnected node or no nodes is called a {\it connected
component}. In this section we derive the distribution of the
number of components in the network.

Let $\{\geq j\}$ denote the network comprising of the ordered nodes
$X_{(j)}, \ldots, X_{(n)}$.  Let $\psi_n(j,k),$ $j=1,\ldots, n,$
$k=1,\ldots, n-j+1,$ denote the probability that in an $n$-node
network there are $k$ components in $\{\geq j\}$, $k=1,\ldots, n-j+1$.
To simplify the notation let $\zeta_{i}(n) := \prob{Y_i \leq r} =
(1-e^{-\lambda(n-i)r})$. The following can be easily verified;
\begin{eqnarray}
  \psi_n(j,n-j+1) \ \ = \ \ \prod_{i=j}^{n-1}(1- \zeta_{i}(n)) ,
  \qquad
  \psi_n(j,1) =  \prod_{i=j}^{n-1}\zeta_{i}(n).
  \label{eq:components-initial-condns}
\end{eqnarray}
Note that $k$ components in $\{ \geq j\}$ can occur in one of two
ways; $k$ components in $\{ \geq (j+1) \}$ and nodes $j$ and $j+1$ are
connected, or $(k-1)$ components in $\{\geq (j+1)\}$ and $j$ not
connected to $j+1$. This leads us to state the following.
\begin{property}
  The probability that there are exactly $k$ components in the graph,
  $\psi_n(1,k)$, is obtained by the recursion
  \begin{equation}
    \label{eq:components-recursion}
    \psi_n(j,k) = \zeta_j(n) \psi_n(j+1,k) + (1-\zeta_j(n)) \psi_n(j+1,k-1).
  \end{equation}
  The initial conditions for the recursion will be given by
  Eqn.~\ref{eq:components-initial-condns}.
  \label{prop:components}
\end{property}

We next investigate the convergence in distribution of the number
of components. From Property~\ref{prop:components} we observe
that as $n \to \infty$, the number of components will essentially
be determined by the last few nodes. To derive the limiting
distribution of the number of components, consider the last node
of the first component. Let $\theta_{n,m}$ denote the probability
that node $m$ is the last node of the first component in an
$n$-node network, $1 \leq m \leq n$.

For any fixed $m$, the probability that the last node of the first
component is the $m\nth$ from the origin goes to 0 as $n \to \infty$,
but for $m=n-s$ we can obtain the following.
\begin{eqnarray}
  \theta_s & := & \lim_{n \to \infty} \theta_{n,n-s} \nonumber \\
  & = & \lim_{n \to \infty} \prod_{i=1}^{m-1} (1 - e^{-r \lambda
  (n-i)}) e^{-r \lambda(n-m)} \nonumber \\
  & = & \lim_{n \to \infty} \frac{P^c_n e^{-r
      \lambda s}}{\prod_{i=1}^{s-1}\ (1-e^{-i r \lambda})}
      \nonumber \\
  & = &  \frac{P_c e^{-r
      \lambda s}}{\prod_{i=1}^{s-1}\ (1-e^{-i r \lambda})},
\end{eqnarray}
where the last equality follows from
Theorem~\ref{thm:conn-asymptote}. As $s \to \infty$, the
denominator decreases monotonically to $P_c$ and $\theta_s$ goes
to zero as $e^{-\lambda r s}$. To obtain the limiting probability
of having $k$ components in the network, conditional on the first
component ending at $m=n-s$, we need $k-1$ components for the
network composed of nodes $n-s+1,\ldots,n$.  The distribution of
the internodal distance between the ordered nodes $n-s+1, \ldots,
n$ is exponential with parameters, $s\lambda,
(s-1)\lambda,\ldots,\lambda$. This is the same internodal
distribution obtained when $s$ nodes are distributed by choosing
their distances from the origin to be exponentially distributed
with mean $1/\lambda$. Thus we can write the following recursive
expression for the limiting probability of the network having $k$
components.
\begin{displaymath}
  \psi(1,k) := \lim_{n \to \infty} \psi_{n}(1,k) = \sum_{s=k}^{\infty}
  \theta_s \psi_s(1,k-1).
\end{displaymath}
$P_c/({\prod_{i=1}^{s-1}\ (1-e^{-i r \lambda})})$ and
$\psi_s(1,k)$ are both bounded sequences. Hence, the series on
the right hand side above converges.  We have thus proved the
following result.
\begin{theorem}
  For fixed $\lambda r$, the number of components in the graph
  converges in distribution, i.e., the probability mass function for
  the number of components in the network converges as $n \rightarrow
  \infty$.
  \label{thm:comp-asymptote}
\end{theorem}
The {\it size} of a component is the number of nodes in that
component. We now derive an expression for the distribution of the
number of components of size $m$. In a network with $n$ nodes,
let $P_m^n(i,k)$ denote the probability that, in $\{\geq i \}$,
there are $k$ components, each of size $m$.  We are interested in
$P_m^n(1,k).$ It is clear that if $mk > n-i+1$, $P_m^n(i,k) = 0$.
Else,
\begin{eqnarray*}
  P_m^n(n-m+1,0) &=& 1-\prob{Y_{n-m+1} \leq r,\ldots, Y_{n-1} \leq r},
  \\ P_m^n(n-m+1,1) &=& \prob{Y_{n-m+1} \leq r,\ldots, Y_{n-1} \leq r}.
\end{eqnarray*}
Conditioning on the location of the first $j \geq i$ such that $Y_j >
r,$ we obtain a recursive relation for $P_m^n(i,k)$ as
\begin{eqnarray}
  P_m^n(i,k) &=& \sum_{j=i+1,j\neq m+i}^{n-km+1} \prob{Y_i \leq r,\ldots,
  Y_{j-2} \leq r, Y_{j-1} > r} P_m^n(j,k) \nonumber \\
  & & \hspace{-15pt} + \ \prob{Y_i
  \leq r,\ldots,Y_{i+m-2} \leq r, Y_{i+m-1} > r} P_m^n(m+i,k-1).
  \label{eq:dist_sizem_comp}
\end{eqnarray}
When $m=1,$ the first factor in the second term above should be
interpreted as $\prob{Y_i > r}.$  The boundary conditions for the
above recursion will be given by
\begin{displaymath}
  P_m^n(i,0) =   \sum_{j=i,j \neq i+m-1}^{n-m}  \prob{Y_i \leq r,\ldots,Y_{j-1}\leq r,Y_j>r}
  P_m^n(j+1,0),
\end{displaymath}
and
\begin{displaymath}
   P_m^n(n-km+1,k)    =  P_m^n(n-(k-1)m+1,k-1) \prob{Y_{n-km+1} \leq
    r,\ldots ,Y_{n-(k-1)m} > r}.
\end{displaymath}

Following the same arguments as in the proof of
Theorem~\ref{thm:comp-asymptote}, we can derive the limiting
distribution of the number of size $m$ components.
\begin{theorem}
  For a fixed $\lambda r$, the limiting distribution of the number of
  size $m$ components is given by the following equation.
  \begin{equation}
    P_m(k) = \lim_{n \rightarrow \infty} P_m^n(1,k) = \sum_{s=mk}^{\infty}
    \theta_s P_m^s(1,k),
    \label{eq:asy_dist_sizem_comp}
  \end{equation}
  where $P_m^s(1,k)$ are as given by (\ref{eq:dist_sizem_comp}).
  \label{thm:dist_size_m_comp}
\end{theorem}

By taking $m=1$ in (\ref{eq:asy_dist_sizem_comp}), we obtain the
asymptotic distribution of the number of isolated nodes in the
network.

\section{Completely Covered Nodes}
\label{sec:redundant}

If there are $k$ nodes in the interval $(X_{(i)}, X_{(i)}+r)$,
then $(k-1)$ are redundant while the $k\nth$ one is necessary for
connectivity and we will say that $k-1$ nodes are `covered' by
node $i$.  From a sensor network perspective, the first $k-1$
nodes in the range of node $i$ to its right may be said to be
redundant. We now determine the distribution of the number of
such covered or redundant nodes in the network. Let $\phi(j,k)$,
$j=1,2,\ldots,n$ and $k=0,1,\ldots,n-j-1,$ denote the probability
that there are $k$ redundant nodes in the network after the
$j\nth$ node, given that the $n$-node network is connected. The
network being connected is denoted by event $C$. We derive a
recursive formula for $\phi(j,k)$ by conditioning on the location
of the last node within the range of the $j-$th node. Our interest
is in $\phi(1,k),$ $k=1,2,\ldots ,n-2$.
\begin{displaymath}
  \phi(j,k)=\sum_{i=j+1}^{j+k+1}
  \prob{X_{(i)} \leq X_{(j)}+r < X_{(i+1)}|C} \phi(i, k-i+j+1),
\end{displaymath}
with boundary condition
\begin{displaymath}
  \phi(j,n-j-1) = \prob{X_{(n)} - X_{(j)} \leq r|C}.
\end{displaymath}
$\prob{X_{(i)} \leq X_{(j)}+r < X_{(i+1)}|C}$ is obtained as follows.
\begin{eqnarray}
  \!\!\!\! \prob{\!X_{(i)}\! \leq \!X_{(j)}\!+\!r \leq X_{(i+1)}\!\mid\!
    C\!}  \!\!\! & \!\!\! = \!\!\! & \!\!\! \prob{\!(X_{(i)}-X_{(j)}\!
    \leq \!r\!) \!\cap \! (X_{(i+1)}-X_{(j)}\!>\!r\!) \!\mid  \!C\!}
  \nonumber\\
  \!\! & \!\! = \!\! & \!\!
  \prob{\!(Y_j+\ldots+Y_{i-1} \!  \leq  \! r \!) \! \cap
    \!(Y_j+\ldots+Y_i \! > \! r \!)  \! \mid \! C \!}\nonumber\\
  \!\! & \!\! = \!\! & \!\! \prob{(Z_{j,i}  \leq  r) \cap (Z_{j,i} +
    Y_i>r ) \mid C} ,
  \label{eq:node between j,i}
\end{eqnarray}
where $Z_{j,i} = Y_j+\ldots+Y_{i-1} $.  Since $Z_{j,i}$ is the sum of
$j-i+1$ exponentials, its density, $g_{Z_{j,i}}(z)$, is given by
\begin{displaymath}
  g_{Z_{j,i}}(z) = \sum_{h=j}^{i-1}\prod_{(m=j,m\neq h)}^{i-1}
  {\frac{n-m}{h-m}}\lambda_h e^{-\lambda_h z},
\end{displaymath}
where $\lambda_h = (n-h) \lambda$ (\cite{Ross01}, Section~5.2.4).
Using this and (\ref{eq:node between j,i}), we get
\[ \prob{X_{(i)} \leq X_{(j)}+r \leq X_{(i+1)}\mid C}  = \hspace{2.5in} \]
\begin{equation}
{\frac{\sum_{h=j}^{i-1}\left( {\frac{n-i}{h-i}}(e^{-\lambda_h
          r}-e^{-\lambda_i r})-e^{-\lambda_h}(1 - e^{-\lambda_i
          r})\right)\prod_{(m=j,m\neq
        h)}^{i-1}{\frac{n-m}{h-m}}}{\prod_{m=j}^i(1-e^{-l(n-i)r})}} .
\end{equation}

Using the initial condition that $\phi(j,k)=0$ for $k > n-j-1$, the
$\phi(j,k)$ can be calculated in the sequence $\phi(n-2,1)$,
$\phi(n-3,1)$, $\phi(n-3,2)$, $\ldots$.

\section{Expected Node Degree}
\label{sec:expected-node-degree}

The degree of a node is the number of nodes lying in its range.
Given a node at $x$, let $p(x)$ denote the probability that
another node is located within distance $r$ of $x.$ While
computing the expected number of nodes of degree $k,$ where $k$
is a fixed integer, we ignore the contribution to the expectation
from nodes lying in $[0,r]$, since as $n \rar \infty$ this
contribution becomes negligible. This will happen since the number
of nodes that fall in $[0,r)$ will approach $\infty$ and thus the
vertex degrees of these nodes for fixed $r$ will tend to
$\infty.$ Let $W_{n,k}$ be the number of nodes of degree $k,$ $k
= 0,1, \ldots$ in an $n-$node network.

\begin{theorem}
  For fixed $\lam r$ as $n \rar \infty,$ $\lim_{n \rar \infty}
  \EXP{W_{n,k}} = c^{-1},$ where $c = (e^{\lam r}-e^{-\lam r})$ and
  the limit is independent of $k.$
  \label{thm:expected-degree-asym}
\end{theorem}

\proof{ Let $X$ be an exponential random variable with parameter
$\lambda.$ Define,
\begin{displaymath}
  p(x) = \prob{x-r \leq X \leq x+r} = ce^{-\lam x}, \qquad x\geq r.
\end{displaymath}
We use the notation $f(n) \sim g(n)$ to indicate that $f(n)/g(n)
\rar 1,$ as $n \rar \infty.$ Since the $n$ nodes are identically
distributed, $\EXP{W_{n,k}}$ will be $n$ times the probability
that any one node in the network has degree $k.$ Condition on this
node being at $x$. Then, the number of nodes lying in $(x-r,x+r)$
is binomial with parameters $(n-1)$ and $p(x).$  By the remark
preceding the statement of the theorem, we ignore the
contribution coming from this node lying in $[0,r).$ Hence,
\begin{eqnarray}
  \EXP{W_{n,k}} & \sim & n {\binom{n-1}{k}} \int_r^{\infty} p(x)^k
  (1-p(x))^{n-k-1} \lam e^{-\lam x} dx \nonumber \\
  & = & n {\binom{n-1}{k}} \int_r^{\infty} c^k
  e^{-\lam k x} (1 - c e^{-\lam x})^{n-k-1} \lam e^{-\lam x} dx \nonumber\\
  & \sim & \frac{n^{k+1}}{k! c} \int_0^{ce^{-\lam r}} y^k
  (1-y)^{n-k-1} dy \nonumber \\
  & = &\frac{n^{k+1}}{k! c} \int_0^{1-e^{-2\lam r}} y^k
  (1-y)^{n-k-1} dy \nonumber \\
  & = &\frac{n^{k+1}}{k! c} \left( \int_0^1 y^k(1-y)^{n-k-1}dy
  - \int_{1-e^{-2\lam r}}^1 y^k(1-y)^{n-k-1} dy \right).
  \label{W_n_k}
\end{eqnarray}
We have used the fact that ${\binom{n-1}{k}} \sim (n-1)^k/k!$ in
deriving the second relation above.  Consider the second integral
in the last equation above. The function $y^k(1-y)^{n-k-1}$ has a
unique maximum in $[0,1]$ at $k/(n-1)$ which tends to zero as $n
\rar \infty.$ Further, the function is monotonically decreasing
in $(k/(n-1),1)$ Thus the second term in the last equation above
is bounded by
\[ {\frac{n^{k+1}}{k! c}} e^{-2\lam r} (1-e^{-2\lam r})^k
(1-e^{-2\lam r})^{n-k-1}, \]
which goes to zero as $n \rar \infty.$ The first term in
(\ref{W_n_k}) is
\[ \frac{n^{k+1}}{k! c} Be(k+1,n-k) = \frac{n^{k+1}}{k! c}
\frac{\Gamma(k+1) \Gamma(n-k)}{\Gamma(n+1)} = \frac{n^{k+1}} {c
n(n-1)\ldots(n-k)},\]
which converges to $c^{-1}$, and hence $\EXP{W_{n,k}} \rar c^{-1}$
as $n \rar \infty.$}

\section{Span and Uncovered Segments}
\label{sec:span-asymptote}

In $G_n(\lambda, r)$, if $Y_i > r$ we can say that there is a portion
between ordered nodes $i$ and $(i+1)$ that is not `covered' and that
there is a hole of size $Y_i - r$. If we think of the nodes as sensors
with range $r$ deployed along a border, then an intruder passing
through the hole will go undetected.  Denoting the length of the hole
between the nodes $i$ and $(i+1)$ by $Z_i$ we have
\begin{displaymath}
  Z_i =\max\{ Y_i - r, 0 \}.
\end{displaymath}
The total length of the holes in the network is then $H(n,r) :=
\sum_{i=1}^n Z_i$ and the number of holes is $\NH(n,r) : = \sum_{i =
  1}^n I_{\{Y_i > r\}}$.

Let $S_n = X_{(n)} - X_{(1)}$ be the span of the network. Since the
$X_i$ are exponentially distributed, as $n \to \infty$, $S_n \to
\infty$ almost surely. However, $H(n,r)$, the total length of the
holes and $\NH(n,r)$, the total number of holes in the network
converge to a proper random variable in distribution.

\begin{theorem}
  As $n \rar \infty,$ $H(n,r)$ and $\NH(n,r)$ converge in distribution
  to random variables with finite mean and variance.
  \label{thm:holes-convergence}
\end{theorem}

\proof{ First, consider the mean and variance of $H(n,r)$ as $n \to
  \infty$.  The density of $Z_i$ is a shifted exponential for $z > 0$
  with a point mass at $0$. Thus the density of $Z_i$, $f_{Z_i}(z)$,
  can be written as
  \begin{displaymath}
    f_{Z_i}(z) \ =  \ \left( 1-e^{-(n-i)\lambda r}\ \right) \
    \delta(z) \ + \ (n-i)\lambda \ e^{-(n-i)\lambda (z+r)}.
  \end{displaymath}
  where $\delta(z)$ is the Dirac-delta function.  The mean and
  variance of $Z_i$ can be shown to be given by
  \begin{eqnarray*}
    \EXP{Z_i} & = & \frac{e^{-(n-i)\lambda r}}{(n-i)\lambda}, \\
    \VAR{Z_i} & = & \frac{e^{-(n-i)\lambda r}(1-e^{-(n-i)\lambda
        r})}{((n-i)\lambda)^2}.
  \end{eqnarray*}
Since $Y_1,\ldots,Y_{n-1}$ are independent, so are the random
variables $Z_1,\ldots ,Z_{n-1}$. The mean and variance of $H(n,r)$
are then given by
  \begin{eqnarray}
    \EXP{H(n,r)} & = & \sum_{k=1}^{n-1}
    \frac{e^{-(n-k)\lambda r}}{(n-k)\lambda} \ = \
    \sum_{k=1}^{n-1} \frac{e^{-k\lambda r}}{k\lambda}
    \label{eq:EXP-H-n-r} , \\
    \VAR{H(n,r)} & = & \sum_{k=1}^{n-1} \frac{e^{-(n-k)\lambda r}
      \left(1-e^{-(n-k)r\lambda}\right)}{((n-k)\lambda)^2} \ = \
    \sum_{k=1}^{n-1} \frac{e^{-k\lambda r} \left( 1-e^{-k\lambda
          r}\right)}{(k\lambda)^2} .
    \label{eq:VAR-H-n-r}
  \end{eqnarray}
  Applying the ratio test to the series in (\ref{eq:EXP-H-n-r}) and
  (\ref{eq:VAR-H-n-r}), we see that $\EXP{H(n,r)}$ and $\VAR{H(n,r)}$
  converge as $n \to \infty$.  Observe that since the variance of
  $H(n,r)$ converges to a finite limit, the usual central limit
  theorem will not be applicable.

  To show convergence in distribution, we must show that the sequence
  of random variables $\{H(n,r)\}$ is tight (relatively compact) and
  that the Laplace transform of $H(n,r)$ converges (see Lemma 2,
  pp.~323 in \cite{Shiryayev84}).  Tightness means that the
  probability of the $H(n,r)$ lying outside a compact set can be made
  arbitrarily small. Tightness also implies that any subsequence
  $H(n_k,r)$ of $H(n,r)$ will contain a subsequence that converges in
  distribution.  We need to show tightness because of the absence of a
  nice closed form expression for the characteristic function of
  $H(n,r)$.  Convergence of the Laplace transform implies uniqueness
  of these limits thereby implying convergence in distribution.

  To show tightness we need to show that for any $\epsilon > 0$, there
  exists a $K>0$ such that $\sup_{n \geq 1} \prob{H(n,r) > K} <
  \epsilon$. $H(n,r)$ are nonnegative random variables and we can use
  Markov inequality to write, for any $K > 0$,
  \begin{displaymath}
    \prob{H(n,r) > K}  \ \le \ \frac{\EXP{H(n,r)}}{K}.
  \end{displaymath}
  Since $\EXP{H(n,r)}$ converges and is finite, for any $\epsilon$, a
  sufficiently large $K$ can be found such that $\EXP{H(n,r)}/K <
  \epsilon$. Thus the random variables $H(n,r)$ are tight.

  To complete the proof of convergence in distribution of $H(n,r)$, we
  have to show that the Laplace transform $L_n(\theta)$ of $H(n,r),$
  converges in some neighborhood of zero.
  \begin{eqnarray*}
    L_n (\theta) & := & \EXP{e^{\theta H(n,r)}} \ = \ \EXP{e^{\theta
        \sum_{i=1}^{n-1} Z_i}} \ = \ \prod_{i=1}^{n-1} \EXP{e^{\theta
        Z_i}}\\ & = & \prod_{i=1}^{n-1} \left( 1 + \frac{\theta
      e^{-(n-i)\lambda r}}{(n-i)\lambda - \theta} \right) \qquad \theta
    < \lambda.
  \end{eqnarray*}
  Taking logarithms on both sides, we get
  \begin{eqnarray*}
    \ln \left( L_n (\theta) \right) & = & \sum_{i = 1}^{n-1} \ln
    \left(1 + \frac{\theta e^{-(n-i) \lambda r}}{((n-i)\lambda -
    \theta)} \right) \\ & = & \sum_{i = 1}^{n-1} \ln \left( 1 +
    \frac{\theta e^{-i r\lambda}}{(\lambda i - \theta)} \right)\\ &
    \leq & \sum_{i = 1}^{n-1} \frac{\theta e^{-i \lambda r}}{(\lambda
    i - \theta)}.
\end{eqnarray*}
The last inequality above is obtained from the inequality
$\ln(1+x) \leq x$. Observe that $\sum_{i = 1}^{\infty} \frac{
\theta e^{-i \lambda
    r}}{(\lambda i - \theta)}$ converges by ratio test.  This proves
the convergence of $L_n(\theta)$ and hence the second part of the
theorem on the convergence of $H(n,r)$ in distribution.

  We now consider convergence in distribution of the number of
  holes. The mean and variance of $\NH(n,r)$ are given by
  {\small
    \begin{eqnarray*}
      \EXP{\NH(n,r)} & = & \sum_{i = 1}^{n-1} e^{- (n-i) \lambda
        r}  \ = \
      \frac{e^{-r\lambda}(1-e^{-(n-1)r\lambda})}{1-e^{-r\lambda}},  \\
      \VAR{\NH(n,r)} & = & \sum_{j = 1}^{n-1}e^{-\lambda_j
        r}(1-e^{-\lambda_j r}) .
    \end{eqnarray*}
  }
  Application of the ratio test shows that both the above series
  converge. Tightness of $\NH(n,r)$ follows by the same argument as
  that used to show the tightness of $H(n,r)$. The Laplace transform
  of $\NH(n,r)$, $J_n(\theta)$, is given by
  \begin{eqnarray}
    J_n(\theta) &  = & \prod_{i=1}^{n-1}\left( 1 - e^{-i \lambda r}
    (1-e^{\theta}) \right) .
  \end{eqnarray}
  Convergence of $J_n(\theta)$ can be shown as for
  $L_n(\theta)$. Thus $\NH(n,r)$ converges in distribution as $n
  \to \infty$.

  This completes the proof of the Theorem. \hfill $\Box$
}

Theorem~\ref{thm:holes-convergence} implies that for large $n$, we can
increase the span of the network over any length with a certain high
probability, by adding more nodes without a corresponding increase in
the length of the holes or the number of holes.

\begin{remark}
  Since the number of components is just one more than the number of
  holes, the convergence in distribution of the number of holes
  follows from Theorem~\ref{thm:comp-asymptote}. Thus this is an
  alternate proof for Theorem~\ref{thm:comp-asymptote}. The limit of
  $J_n(\theta$) can be used to obtain the asymptotic moments for the
  number of components.
\end{remark}

The asymptotic distribution of the span is also known. From
Examples~3.3 and 3.5 of \cite{Castillo87}, we have that $ \lam^{-1}
X_{(1)} \log(n/(n-1))$ converges in distribution to a Weibull
distribution and $\lam X_{(n)} - \log(n)$ converges in distribution to
a Gumbel distribution. This allows us to state the following result
for the asymptotic distribution of the span.
\begin{theorem}
  $\lambda(X_{(n)} - X_{(1)}) - \log(n)$ converges in distribution to
  a Gumbel distribution.
  \label{thm:span-distribution-asym}
\end{theorem}

Thus, the $100(1-\al)\%$ confidence interval for the span based on the
asymptotic distribution will be of the form $\lambda^{-1} (\log(n) \pm
c(\al))$ where $c(\al)$, is independent of $n.$

\section{Strong Law Results}
\label{sec:strong-laws}
In this section we derive almost sure convergence results for the
connectivity and the largest nearest neighbor distances and a
limiting result for the almost surely connected part of the
exponential random geometric graph.

Define $c_n$ and $d_n$ the connectivity and largest nearest neighbor
distances respectively as
\begin{eqnarray}
  c_n &=& \inf \{ r > 0: G_n(\lam,r) \mbox{ is connected} \},
  \label{def_con_dist} \\
  d_n &:= & \max_{1 \leq i \leq n} \min_{1 \leq j \leq n, j \ne i} \{ |X_i - X_j| \}
  \label{def_lnnd}
\end{eqnarray}
\begin{theorem}
  For fixed $\lam > 0,$
  \begin{enumerate}
  \item
    \begin{equation}
      \limsup_{n \rar \infty} \frac{\lam c_n}{\ln (n)} = \limsup_{n \rar
        \infty} \frac{\lam d_n}{ \ln (n)} = 1, \qquad \mbox{almost
        surely.}
      \label{eqn:slln_c_n-d_n}
    \end{equation}
  \item
    \begin{equation}
      \liminf_{n \rar \infty} \frac{\lam \ln(n) c_n}{c} \geq 1,
      \qquad \liminf_{n \rar \infty} \frac{\lam \ln(n) d_n}{c} \geq 1,
      \qquad \mbox{almost surely.}
      \label{eqn:liminf_lower_bound}
    \end{equation}
    where $c = \sum_{j=1}^{\infty} j^{-2}.$

  \item Let $r$ be fixed, $k_n = \lfloor n ( 1- a \ln(n)/n) \rfloor$
    where $\lfloor \cdot \rfloor$ denotes the integer part and $a >
    (\lam r)^{-1}$. Let $G_n(k_n,\lam,r)$ denote the graph
    $G_n(\lam,r)$ restricted to the first $k_n$ ordered points. Then,
    \begin{equation}
      \prob{G_n(k_n,\lam,r) \mbox{ is disconnected infinitely often }} = 0.
      \label{eqn:a.s.connected-part}
    \end{equation}
  \end{enumerate}
  \label{thm:slln_exp}
\end{theorem}

\proof{
  \[
  \prob{c_n \geq y}  =  \prob{\cup_{i=1}^{n-1}\{Y_i \geq y\}}
  \leq \sum_{i=1}^{n-1} e^{-\lam i y} = e^{-\lam y}
  \frac{1-e^{-(n-1)\lam y}} {1-e^{-\lam y}}.
  \]
  Taking $y=(1+\ep)\log(n)/\lam,$ and applying the ratio test, we see
  that
  \begin{displaymath}
    \sum_{n=2}^{\infty} \prob{\lam c_n \geq (1+\ep)\log(n)} < \infty.
  \end{displaymath}
  By the Borel-Cantelli lemma, $\prob{\lam c_n \geq (1+\ep) \log(n)
    \mbox{ i.o. }} = 0.$ Since $\ep > 0$ is arbitrary, we conclude that
  $\limsup (\lam c_n/\log(n)) \leq 1$ a.s.

  To show that the $\limsup$ is exactly equal to one, consider the
  record values denoted by $R_n$ defined as follows: Let $N(1) = 1.$
  For $n \geq 2$, define $N(n) = \inf \{ k > N(n-1): X_k > X_{N(n-1)}
  \}$. Define $R_n = X_{N(n)}.$ Since the exponential density has
  unbounded support, there will be a.s. infinitely many record
  values. Consider the sequence $Z_n = R_n - R_{n-1}$, $n \geq 2$. By
  the memoryless property of the exponential, $Z_n$ is a sequence of
  independent exponential random variables with mean
  $\lam^{-1}$. Since, for any $\ep >0,$
  \[
  \sum _{n=2}^{\infty} \prob{\lam Z_n > (1-\ep)\log(n)} =
  \sum _{n=2}^{\infty} n^{-(1-\ep)} = \infty,
  \]
  it follows from the Borel-Cantelli Lemma that $\limsup \lam
  Z_n/\log(n) = 1$ a.s.  The above result implies that $\limsup (\lam
  d_n/\log(n)) \geq 1$ a.s. by considering the sequence of graphs
  $G_{N(k)}(\lambda, r).$ Part~1 of the theorem now follows because
  $d_n \leq c_n.$

  To prove part~2 for $c_n$, we consider the asymptotic behavior of the
  probability that $G_n(\lambda, r_n)$ is connected for the sequence
  of cutoff distances $r_n = c/(\lam (1+ \ep) \ln(n))$, where $c$ is
  as defined in the theorem statement.
  \[
  P_c^n = \prob{G_n(\lam, r_n) \mbox{ is connected }} =
  \prod_{i=1}^{n-1}(1-\exp(-\lam i r_n)).
  \]
  Taking logarithms and expanding the logarithm, we get
  \begin{eqnarray*}
    \ln(P_c^n) & \ = \ &- \sum_{i=1}^{n-1} \sum_{j=1}^{\infty}
    \frac{e^{-\lam i j r_n}}{j} \ \  = \ \ - \sum_{j=1}^{\infty}
    \frac{e^{-\lam j r_n}(1-e^{-\lam j (n-1)r_n})}{j (1-e^{-\lam j
        r_n})}.
  \end{eqnarray*}
Since $r_n \rar 0,$ and $n\:r_n \rar \infty,$ we have $e^{-\lam j
r_n} \rar 1,$ $1-e^{-\lam j (n-1)r_n} \rar 1,$ and $1-e^{-\lam j
r_n} \sim \lam j r_n.$ Hence,
\[
\ln(P_c^n) \sim  - \frac{1}{\lam r_n} \sum_{j=1}^{\infty}
    j^{-2}.
\]
  Plugging in the expression for $r_n$ we get $P_c^n \sim
  n^{-(1+\ep)},$ which is summable. The result for $c_n$ in Part~2
  now follows from the Borel-Cantelli lemma. To prove part~2 for
  $d_n$ let $y_n = \frac{c}{\lam (1+\ep)\log n},$ and consider,
\begin{eqnarray*}
    \lefteqn{ \prob{ d_n \leq y_n} }\\
    & =  & \prob{\cap_{i=2}^{n-1}((Y_{i-1}
      \leq y_n)\cup (Y_{i} \leq y_n)) \cap (Y_{1} \leq y_n) \cap
      (Y_{n-1} \leq y_n)} \\
    & \leq & \prob{\cap_{i=1}^{\lfloor n/2 \rfloor} ((Y_{2i-1} \leq
      y_n) \cup (Y_{2i} \leq y_n))} \\
    & \leq & \prod_{i=1}^{\lfloor n/2 \rfloor}(1-e^{-\lam (2n - 4i
      -1)y_n}).
  \end{eqnarray*}
Take logarithms on both sides and using the Taylor expansion, we
get,
\begin{eqnarray*}
\ln(\prob{ d_n \leq y_n})
 & = & -\sum_{j=0}^{\infty} \sum_{i=1}^{\lfloor n/2\rfloor}
  \frac{e^{-\lam(2n-4i+1)jy_n}}{j}\\
 & = & -\sum_{j=0}^{\infty} \frac{e^{-\lam jy_n(2n+1)}}{j} \sum_{i=1}^{\lfloor n/2\rfloor}
  (e^{4\lam jy_n})^i \\
 & = & -\sum_{j=0}^{\infty} \frac{1}{j} e^{-\lam j
 y_n (2n -3)} \frac{1 - e^{4 \lam j y_n \lfloor n/2\rfloor}}{1 -
 e^{4 \lam j y_n}} \\
 & = & -\sum_{j=0}^{\infty} \frac{1}{j}\frac{e^{-\lam j
 y_n (2n -3)} - e^{4 \lam j y_n (2n-3-4\lfloor n/2\rfloor)}}{1 -
 e^{4 \lam j y_n}} \\
 & \sim & \frac{-1}{4\lam y_n} \sum_{j=0}^{\infty}
  \frac{1}{j^2},
\end{eqnarray*}
where the last approximation follows since $y_n \rar 0$ and $n
y_n \rar \infty$ which implies that $\exp(-\lam j
 y_n (2n -3)) \rar 0$, $\exp(4 \lam j y_n (2n-3-4\lfloor n/2\rfloor))
\sim \exp(- 12 \lam j y_n) \rar 1,$ and $1 -
 \exp(4 \lam j y_n) \sim - 4 \lam j y_n$. Substituting for $y_n$ we
get,
\begin{equation}
P[d_n \leq \frac{c}{\lam (1+\ep)\log n}] \sim \frac{1}{n^{1+\ep}},
\end{equation}
which is summable. Part~2 of the theorem for $d_n$ now follows
from the Borel-Cantelli Lemma.

To prove part~3, consider %
  \begin{eqnarray*}
    \prob{G_n(k_n,\lam,r) \mbox{ is not connected } } & \ \leq \ &
    \sum_{i=1}^{k_n-1} e^{-\lam r (n-i)} \\
    & \ = \ & \frac{e^{\lam r}}{e^{\lam r}-1} \left(
    e^{-\lam r (n-k_n)} - e^{-\lam r n} \right).
  \end{eqnarray*}

  For large $n$, $n-k_n \sim a \log(n),$ and hence the above
  probability is summable. The result now follows from the
  Borel-Cantelli Lemma.  }

\section{Truncated Exponential Graph}
\label{sec:truncated-exponential}
We now consider the RGG $G_n(\lam,r,T)$ where the nodes are
distributed independently according to the density function
\[
g_{\lam,T}(x) = \frac{\lam e^{-\lam x}}{1 - e^{-\lam T}}, \qquad 0 \leq
x \leq T,
\]
and have a cutoff $r$. This distribution allows us to consider RGGs
with finite support where the distribution of the nodes is not
uniform. We derive asymptotic results for the connectivity and largest
nearest neighbor distances for this graph. \cite{Penrose03} derives
similar results for dimensions $d \geq 2$ and for general densities
having bounded support. \cite{Appel97b,Appel02} show such strong laws
for the uniform RGG for $d \geq 1.$ \cite{Iyer05a} obtains strong law
results for the one dimensional uniform RGG using the graph with
independent exponential spacings of \cite{Mccolm04}.

In deriving our results, note that unlike in the exponential RGG,
the spacings in $G_n(\lam,r,T)$ are not independent. Our proof
technique is as follows. We show that the graph $G_n(\lam,r,T)$
has the same asymptotic behavior as that of a graph $G_n^*$ which
is constructed by considering the first $n$ nodes of an
exponential RGG on $N$ vertices. Here $N = N(n) := \lfloor n/p
\rfloor$ and $p = 1 - \exp(- \lam T)$. Spacings in the graph
$G_n^*$ are independent and hence it is possible to derive
results easily for this graph. This technique allows us to think
about $G_n(\lam,r,T)$ in terms of the graph $G_n^*$ whose
properties can be more easily visualized. This is similar to the
approach of \cite{Mccolm04} for the uniform RGG.

Let $X_1,X_2, \ldots $ be a sequence of independent random variables
with density $g_{\lam,T}$. The vertex set of $G_n(\lam,r,T)$ is $V_n =
\{ X_1,\ldots ,X_n \}.$ Let $N(n)$ be as defined above and let
$Z_1,Z_2,\ldots ,$ be a sequence of exponential random variables with
mean $\lam^{-1}.$ Let $Z_{1,N},\ldots ,Z_{N,N}$ denote the ordered
values of the first $N(n)$ random variables $Z_1,\ldots , Z_N.$ Define
the graph $G_n^*(\lam,r)$ to be the RGG with cutoff $r$ and vertex set
$V_n^* = \{ Z_{1,N}, \ldots , Z_{n,N} \}.$ We denote by
$G_n^*(\lam,r,t)$ the graph with vertex set $V_n^*$ conditioned on
$Z_{n+1,N} = t.$ It is easy to see that the conditional density of
first $n$ ordered observations $Z_{1,N},\ldots ,Z_{n,N}$ given
$Z_{n+1,N}=t$ is given by (see \cite{Rohatgi00}, pp. 175--176),
\begin{equation}
f_{Z_{1,N},\ldots,Z_{n,N}\mid Z_{n+1,N}}(z_1,\ldots,z_n \mid t) =
{\frac{n! \lam^n }{(1-e^{-\lam t})^n}}e^{-\lam
\sum_{i=1}^{n}z_{i}},
\label{eqn:conditional_density_truncated_exp_graph}
\end{equation}
for $0 < z_1 < \ldots < z_ n < t.$ The key observation is that
the above function is also the joint density function of $n$
i.i.d. ordered observations from $g_{\lam,t}$. Further, we have
the following lemma which states that $Z_{n+1,N}$ is close to $T$
with large probability as $n \rar \infty.$ Subsequent to this
lemma, we show that the graphs $G_n$ and $G_n^*$ have the same
asymptotic behavior.

\begin{lemma}
  $Z_{n+1,N} \rar T$ in probability as $n \rar \infty.$
  \label{convergence_to_T}
\end{lemma}

\proof{ We show that the mean and variance of $Z_{n+1,N}$ converge to
  $T$ and $0$ respectively. The result then follows from Chebyshev's
  inequality.
\[ \EXP{Z_{n+1,N}} = \sum_{i=0}^{n} \frac{1}{(N-i) \lam} =
\frac{1}{\lam} \sum_{i=N-n}^N \frac{1}{i}. \]
Hence,
\[ \int_{N-n}^{N+1} \frac{1}{x} dx \  \leq \ \lam \EXP{Z_{n+1,N}}
\ \leq \ \int_{N-n-1}^{N} \frac{1}{x} dx. \]
Both the integrals above converge to $\lambda T$ as $n \rar
\infty$ by the definition of $N.$
\begin{eqnarray*}
  \VAR{Z_{n+1,N}} & \ = \ & \sum_{i=0}^{n} \frac{1}{(N-i)^2 \lam^2} \ = \
  \frac{1}{\lam^2} \sum_{i=N-n}^N \frac{1}{i^2} \\
  & \leq & \frac{1}{\lam^2} \int_{N-n-1}^N \frac{1}{x^2} dx \rar 0.
\end{eqnarray*}
Thus for any $\ep > 0,$ and $n$ sufficiently large, we have $\mid
\EXP{Z_{n+1,N}} - T \mid < \ep/2.$ Hence
\begin{eqnarray*}
  \prob{ \mid Z_{n,N+1} - T \mid > \ep} & \leq & \prob{ \mid
    Z_{n+1,N} - \EXP{Z_{n+1,N}} \mid > \ep/2} \\
  & \leq & \frac{4\VAR{Z_{n+1,N}}}{\ep^2} \rar 0.
\end{eqnarray*}
This completes the proof the lemma. \qed
}

We now show that the graphs $G_n$ and $G_n^*$ have the same
thresholding behavior. To do this we need some notations.

\begin{definition}
  If $A$ and $B$ are graphs such that $A$ and $B$ share the same
  vertices, and the edge set of $A$ is a subset of the edge set of
  $B$, we will write $A \leq B.$ Let $\Theta$ be a property of a
  random geometric graphs such that if $A \leq B$ and $A \in \Theta,$
  then $B \in \Theta.$ (Here $A \in \Theta$ is used to denote that RGG
  $A$ has property $\Theta$.) Then $\Theta$ is called an
  ``upwards-closed'' property. If $B \in \Theta$ implies $A \in \Th$, then
  $\Theta$ is said to be a ``downwards-closed'' property.
\end{definition}

Fix an upwards-closed property $\Theta.$ For any two functions
$\delta, \gamma: Z^+ \rar \Re^+$, we write $\delta \ll \gamma$
(resp. $\del \gg \gamma$) if $\del(n)/\gamma(n) \rar 0,$ (resp.
$\gamma(n)/\delta(n) \rar 0$) as $n \rar \infty.$ In what follows
we will write $\delta$ for $\delta(n).$ Let $G_n(r)$, be any
random geometric graph on $n$ vertices with cutoff $r$.

\begin{definition}
  A function $\delta_{\Theta}: Z^+ \rar \Re^+$ is a weak threshold
  function for $\Theta$ if the following is true for every function
  $\delta: Z^+ \rar \Re^+,$
  \begin{itemize}
  \item if $\delta(n) \ll \delta_{\Theta}(n)$, then $\prob{G_n(\delta)
    \in \Theta} = o(1),$ and
  \item if $\del(n) \gg \del_{\Theta}(n)$ then $\prob{G_n(\del) \in
    \Theta} = 1 - o(1).$
  \end{itemize}

  A function $\del_{\Theta}: Z^+ \rar \Re^+$ is a strong threshold
  function for $\Theta$ if the following is true for every fixed
  $\ep>0,$
  \begin{itemize}
  \item if $\prob{G_n((1-\ep)\del_{\Theta}) \in \Theta} = o(1),$ and
  \item if $\prob{G_n((1+\ep)\del_{\Theta}) \in \Theta} = 1 - o(1).$
  \end{itemize}
\end{definition}

Before proceeding further, we show the following two monotonicity
properties that will be used subsequently.

\begin{lemma} Let $\Th$ be any upwards-closed property. Then,
  \begin{enumerate}
  \item For any $0 < T_1 < T_2$, $\prob{G_n(\lam,\del,T_1) \in \Th}
    \geq \prob{G_n(\lam,\del,T_2) \in \Th}.$
  \item For any $0 < \lam_1 < \lam_2$, $\prob{G_n^*(\lam_1,\del) \in
    \Th} \leq \prob{G_n^*(\lam_2,\del) \in \Th}.$
  \item For any $c > 0,$ $\prob{G_n(\lam,\del,T) \in \Th} =
    \prob{G_n(c^{-1}\lam,c\del,cT) \in \Th}$.
  \end{enumerate}
  \label{lemma_monotonicity}
\end{lemma}
\proof{

  Let $U_{(1)},U_{(2)},\ldots ,U_{(n)}$ be $n$ ordered uniform random
  variables on $(0,1).$ The ordered vertex sets of the graphs
  $G_n(\lam,\del,T_1)$, $G_n(\lam,\del,T_2)$, $G_n^*(\lam_1,\del)$ and
  $G_n(\lam_2,\del)$ may be defined using the ordered uniform variables as follows:
  $V_1 = \{ - \frac{1}{\lam} \ln (1 - U_{(i)}(1-e^{-\lam T_1})) \}_{i=1}^n,$
  $V_2 = \{ - \frac{1}{\lam} \ln (1 - U_{(i)}(1-e^{-\lam T_2}))
  \}_{i=1}^n,$ $V_3 = \{ - \frac{1}{\lam_1} \ln(1-U_{(i)})
  \}_{i=1}^n,$ and $V_4 = \{ - \frac{1}{\lam_2} \ln(1-U_{(i)})
  \}_{i=1}^n.$ If we denote the respective edge sets by $E_i,$
  $i=1,\ldots ,4,$ then it is easy to see that $E_2 \subset E_1$ and
  $E_3 \subset E_4.$ The result for parts~1~and~2 now follows from the
  definition of an upwards-closed property. To prove part~3,
  observe that $\Th$ being an upwards-closed property depends
  only on the existence of edges between certain pairs of nodes in
  a given configuration of vertices $(Z_1,\ldots ,Z_n)$ of the graph
  $G_n(\lam,\del,T).$ Let $A \subset \Re^n$ be such that $\Th$ holds
  whenever $(Z_1,\ldots ,Z_n) \in A$ at
  cutoff $\del$, then clearly, it holds at cutoff $c \del$ if
  $(Z_1,\ldots ,Z_n) \in cA$. The joint density of $n$ independent truncated
  exponential random variables on $[0,T]$ is given by
  (\ref{eqn:conditional_density_truncated_exp_graph}), with $t$ replaced by
  $T.$ Hence,
\begin{equation}
\prob{G_n(\lam,\del,T) \in \Th} = \int_{\{ (z_1,\ldots ,z_n) \in
A \}} \frac{n! \lam^n }{(1-e^{-\lam T})^n}e^{-\lam
\sum_{i=1}^{n}z_{i}} dz_1 \ldots dz_n. \label{prob_g_n_del_theta}
\end{equation}
From the above remarks on upwards-closed property and
(\ref{prob_g_n_del_theta}), we have
\begin{eqnarray*}
\prob{G_n(c^{-1}\lam,c\del,cT) \in \Th} & = & \int_{ \{
(z_1,\ldots ,z_n) \in cA \} } \frac{n! c^{-n}\lam^n
}{(1-e^{-c^{-1}\lam cT})^n}e^{-c^{-1}\lam \sum_{i=1}^{n}z_{i}} dz_1 \ldots dz_n \\
 & = & \int_{\{ c^{-1}(z_1,\ldots
,z_n) \in A \}} \frac{n! c^{-n}\lam^n }{(1-e^{-\lam
T})^n}e^{-c^{-1}\lam \sum_{i=1}^{n}z_{i}} dz_1 \ldots dz_n \\
 & & \\
 & & \qquad \mbox{ Change variables } cu_i =
z_i, \;\; i=1,\ldots , n. \\
 & & \\
 & = & \int_{\{ (u_1,\ldots ,u_n) \in
A \}} \frac{n! \lam^n }{(1-e^{-\lam T})^n}e^{-\lam
\sum_{i=1}^{n}u_{i}} du_1 \ldots du_n  \\
 & = & \prob{G_n(\lam,\del,T) \in \Th},
\end{eqnarray*}
which proves part~3. }

\begin{lemma}
  Let $\del : Z^+ \rar \Re^+$, $T > 0$ and $\al \in
  (0,1).$ Let $G_n(\lam,\del,T)$ and $G_n^*(\lam,\del)$ be the
  random geometric graphs defined above. Then for all $n$
  sufficiently large the following hold.
  \begin{enumerate}
  \item If $\prob{G^*_n(\lam,\del)\in \Theta} \leq \alpha,$ then
    $\prob{G_n(\lam,(1-\ep)\del,T) \in \Theta} \leq
    {\frac{\alpha}{1-\alpha}}.$
  \item If $\prob{G_n(\lam,\del,T) \in \Theta} \leq \alpha,$ then
    $\prob{G^*_n(\lam,(1-\ep)\del) \in \Theta}\leq 2\alpha.$
  \item If $\prob{G^*_n(\lam,\del)\in \Theta} \geq 1-\alpha,$ then
    $\prob{G_n(\lam,(1+ \frac{\ep}{1-\ep})\del,T) \in \Theta} \geq
    {\frac{1-2\alpha}{1-\alpha}}.$
  \item If $\prob{G_n(\lam,\del,T) \in \Theta} \geq 1- \alpha,$ then
    $\prob{G^{*}_n(\lam,(1+\ep)\del) \in \Theta} \geq 1- 2\alpha.$
  \end{enumerate}
  \label{lemma_threshold_truncated_exp}
\end{lemma}
\proof{

  Let $Z_{(n+1)} = Z_{n+1,N}$, be the random variables defined prior to
  Lemma~\ref{convergence_to_T}. For any $\epsilon,\alpha \geq 0,$
  from Lemma~\ref{convergence_to_T},
  there exists a $M \geq  0$ such that, for all $n \geq M,$
  \begin{equation}
    \prob{\mid Z_{(n+1)} - T \mid \leq \ep} \geq 1-\alpha
    \label{interval for Z_n+1}
  \end{equation}
  For the sake of simplicity we will take $T=1$ in this proof and
  write $G_n(\lam,\del)$ for $G_n(\lam,\del,1).$
  \begin{eqnarray*}
    \lefteqn{\prob{G_n^*((1-\ep)\lam,\del) \in \Theta} } \\
    & = & \int_0^{\infty} \prob{ G_n^*((1-\ep)\lam,\del) \in \Theta \mid
      Z_{(n+1)}=z} f_{Z_{(n+1)}}(z)dz\\
    & \geq & \int_{1-\ep}^{1+\ep} \prob{G_n^*((1-\ep)\lam,\del) \in
      \Theta \mid Z_{(n+1)}=z} f_{Z_{(n+1)}}(z)dz\\
    & = & \int_{1-\ep}^{1+\ep} \prob{G_n((1-\ep)\lam,\del,z) \in
      \Theta} f_{Z_{(n+1)}}(z)dz\\
    & \geq & \prob{G_n((1-\ep)\lam,\del,(1+\ep)) \in \Theta}
    \int_{1-\ep}^{1+\ep} f_{Z_{(n+1)}}(z) dz.
  \end{eqnarray*}
  In deriving the second inequality above we have used
  Lemma~\ref{lemma_monotonicity}(1). Using (\ref{interval
    for Z_n+1}), we get
  \[
  \prob{G_n^*((1-\ep)\lam,\del) \in \Theta} \geq (1-\alpha)
  \prob{G_n((1-\ep)\lam,\del,(1+\ep)) \in \Theta},
  \]
  for all $n \geq M$.  Since $\prob{G^*_n(\lam,\del) \in \Theta} \leq
  \alpha,$ it follows from Lemma~\ref{lemma_monotonicity}(2) that

  \begin{eqnarray*}
    \prob{G_n^*((1-\ep)\lam,\del) \in \Theta} & \leq &  \alpha, \\
    \prob{G_n((1-\ep)\lam,\del,(1+\ep)) \in \Theta} & \leq &
    \frac{\alpha}{1-\alpha} .
  \end{eqnarray*}
  Hence, by Lemma~\ref{lemma_monotonicity}(3),
  \[
  \prob{G_n(\lam,(1-\ep)\del, (1+\ep)(1-\ep)) \in
    \Theta } \leq {\frac{\alpha}{1-\alpha}} .
  \]
  Since $(1+\ep)(1-\ep) \leq 1$, part~1 of the lemma follows from another
  application of Lemma~\ref{lemma_monotonicity}(1).

  To prove part~2, consider
  \begin{eqnarray*}
    \lefteqn{\prob{G_n^*((1-\ep)^{-1}\lam,(1-\ep)\del)\in \Theta} } \\
    & = & \int_0^{\infty} \prob{ G_n^*((1-\ep)^{-1}\lam,(1-\ep)\del)\in
      \Theta \mid Z_{(n+1)}=z} f_{Z_{(n+1)}}(z) dz \\
    & \leq & \alpha + \int_{1-\ep}^{1+\ep}
    \prob{G_n^*((1-\ep)^{-1}\lam,(1-\ep)\del) \in \Theta \mid
      Z_{(n+1)}=z} f_{Z_{(n+1)}}(z) dz \\
    & = & \alpha + \int_{1-\ep}^{1+\ep}
    \prob{G_n((1-\ep)^{-1}\lam,(1-\ep)\del,z) \in \Theta}
    f_{Z_{(n+1)}}(z) dz \\
    & \leq & \alpha + \prob{G_n((1-\ep)^{-1}\lam,(1-\ep)\del,(1-\ep)) \in
      \Theta}
    \int_{1-\ep}^{1+\ep} f_{Z_{(n+1)}}(z)dz, \\
    & \leq & \alpha + \prob{G_n(\lam,\del,1) \in \Theta},
  \end{eqnarray*}
  for all $n \geq N$. The first inequality above follows from
  (\ref{interval for Z_n+1}). The second inequality follows from
  Lemma~\ref{lemma_monotonicity}(1) while the last inequality follows
  from Lemma~\ref{lemma_monotonicity}(3). From the given condition
  $\prob{G_n(\lam,\del,1) \in \Theta} \leq \alpha$, it follows that %
  \begin{eqnarray*}
    \prob{ G_n^*((1-\ep)^{-1} \lam,(1-\ep)\del) \in \Theta} & \leq &
    \alpha + \alpha = 2\alpha.
  \end{eqnarray*}
  Part~2 now follows from Lemma~\ref{lemma_monotonicity}(2).

  To prove part~3, proceeding as above, we get
  \begin{eqnarray*}
    \lefteqn{\prob{G_n^*((1-\ep)^{-1}\lam,\del) \notin \Theta} } \\
    & = & \int_0^{\infty} \prob{G_n^*((1-\ep)^{-1}\lam,\del)
      \notin \Theta \mid Z_{n+1}=z} f_{Z_{(n+1)}}(z)dz  \\
    & \geq & \int_{1-\ep}^{1+\ep} \prob{G_n^*((1-\ep)^{-1}\lam,\del)
      \notin \Theta \mid Z_{n+1}=z} f_{Z_{(n+1)}}(z) dz \\
    & = & \int_{1-\ep}^{1+\ep} \prob{G_n((1-\ep)^{-1}\lam,\del,z)
      \notin \Theta}
    f_{Z_{(n+1)}}(z) dz \\
    & \geq & \prob{G_n((1-\ep)^{-1}\lam,\del,(1-\ep)) \notin \Theta}
    \int_{1-\ep}^{1+\ep}f_{Z_{(n+1)}}(z) dz, \\
    & \geq & (1-\alpha) \prob{G_n((1-\ep)^{-1}\lam,\del,1-\ep) \notin \Theta},
  \end{eqnarray*}
  for all $n \geq M.$  Since the given condition
  $\prob{G_n^*(\lam,\del) \notin \Theta} \leq \alpha,$ implies \\
  $\prob{G_n^*((1-\ep)^{-1}\lam,\del) \notin \Theta} \leq \alpha,$ it
  follows that
  \[
  \prob{G_n((1-\ep)^{-1}\lam,\del,1-\ep) \notin \Theta} \leq
    {\frac{\alpha}{1-\alpha}}
  \]
  This implies that
  \[
  \prob{G_n(\lam,(1 + \frac{\ep}{1-\ep})\del) \in \Theta } \geq
  \frac{1-2\alpha}{1-\alpha} ,
  \]
  and we have the proof for part~3.

  To prove part~4 we proceed as above to get the following inequality.
  \[
  \prob{G_n^*(\lam,(1+\ep)\del) \notin \Th }  \leq  \al +
  \prob{G_n(\lam, (1+\ep)\del,1+\ep) \notin \Th} \leq  2 \al.
  \]

  This completes the proof the lemma.
}

\begin{remark}
  The above results extend to downwards-closed properties as well.
\end{remark}

The following theorem is now an easy corollary of the above Lemma.

\begin{theorem}
  The sequence of random geometric graphs $G_n(\lam,\del,T)$ and
  $G_n^*(\lam,\del)$ have the same weak and strong thresholds.
  \label{thm:threshold_equivalence_truncated_exp}
\end{theorem}

\begin{remark}
  Proceeding as in the proof of
  Lemma~\ref{lemma_threshold_truncated_exp}, we can show that %
  \[
  \mid \prob{G_n(\lam,\del,T) \in \Th} - \prob{G_n^*(\lam,\del)
    \in \Th} \mid \rar 0, \qquad n \rar \infty.
  \]
  This implies that \emph{all} the asymptotic probabilities of
  $G_n(\lambda, \delta, T)$ satisfying any monotone property (upward or
  downward closed) can be obtained by studying the corresponding
  probabilities for $G_n^*(\lambda,\delta).$
\end{remark}

We now use Theorem~\ref{thm:threshold_equivalence_truncated_exp} to
derive the threshold probability for connectivity and strong law
results for the connectivity and largest nearest neighbor distance for
the graph $G_n.$ We will, without any further reference to the above
theorem, work with the graph $G_n^*$ instead of $G_n$.

\begin{theorem} Let $p = 1 - \exp(-\lam T).$ Then the sequence of edge
distances $\del(n) = \frac{p}{\lam(1-p)} \frac{\ln(n)}{n}$ is a
strong (and weak) threshold for connectivity for the graph
$G_n(\lam,\del,T).$
\label{thm:threshold_connectivity_truncated_exp}
\end{theorem}

\proof{

  Let $r_n = a \del(n),$ where $a \geq 0$ is a constant. Note that $r_n
  \rar 0$ while $n r_n \rar \infty$ and $n \sim Np.$ Let
  $P_n^{*c}$ be the probability that $G_n^*(\lam,r_n)$ is
  connected. Then,
\[   \ln(P_n^{*c})  =  \sum_{j=N-n}^N \ln(1 - e^{-\lam r_n j} ). \]
Since $j r_n \rar \infty$ for all $j = N-n,\ldots , N,$ using
$\ln(1-x) \sim x$ as $x \rar 0,$ and summing the resultant
geometric series, we get
\[   \ln(P_n^{*c})  \sim  - e^{-\lam r_n (N-n)} \frac{1 - e^{- \lam r_n (n+1)}}{1 -
e^{-\lam r_n}}. \]
Substituting for $r_n = a \del(n)$ while noting that $(N-n)/n
\sim (1-p)/p,$ $1 - e^{- \lam r_n (n+1)} \rar 1$ and $1 -
e^{-\lam r_n} \sim - \lam r_n,$ we obtain,
\[ \ln(P_n^{*c})  \sim - \frac{n^{1-a}}{\ln(n)}. \]
  Thus, for $a=1+\ep$, $P_n^{*c}$ converges to $1$ and converges to
  $0$ for $a = 1-\ep.$ This shows that $\del(n)$ is a strong
  threshold for connectivity for $G_n.$ Similarly one can show that
  $\del(n)$ is a weak threshold as well.
}

\begin{remark}
  Note that $\frac{n \del(n)}{\ln(n)} = \frac{p}{\lam (1-p)}$, where
  $\frac{p}{\lam (1-p)}$ is the reciprocal of the minimum
  $g_{\lam,T}(x).$ Thus the behavior of the distance required to
  connect the graph is determined by the minimum of the density since
  in the vicinity of this point vertices are more sparsely
  distributed. The normalization $\frac{n}{\ln(n)}$ is the same as in
  the case of uniform distribution of nodes.
\end{remark}

We now state a strong law result for the connectivity distance
$(c_n(\lambda,T))$ and the largest nearest neighbor distance
$(d_n(\lambda,T))$ for $G_n(\lambda,\cdot,T)$.  In the following, we
drop the reference to parameters $\lambda$ and $T$ when referring to
$c_n(\lambda,T)$ and $d_n(\lambda,T)$.

\begin{theorem}
  Let $\lam, T > 0.$ The connectivity and largest
  nearest neighbor distances of the graph $G_n(\lam,\cdot,T)$
  satisfy
  \begin{enumerate}
  \item $\lim_{n \rar \infty} \frac{n c_n}{\ln(n)} = \frac{p}{\lam
    (1-p)}$ almost surely.
  \item $\lim_{n \rar \infty} \frac{n d_n}{\ln(n)} = \frac{p}{2\lam
    (1-p)}$ almost surely.
  \end{enumerate}
  \label{thm:slln_trunc_exp}
\end{theorem}

\proof{

  Let $c_n^*$ and $d_n^*$ be the connectivity and largest nearest
  neighbor distances respectively of $G_n^*(\lam,\del).$ Let $Y_{i,n}$
  be the spacings between the vertices in $G_n^*$. The $Y_{i,n}$ are
  independent and exponentially distributed with mean $(\lam
  (N-i))^{-1},$ where $N = \lfloor n/p \rfloor,$ and $p = 1 -
  \exp(-\lam T).$ Note that $(1-p)/p \sim N/n -1.$ Let $y = y_n :=
  (p(1+\ep)\ln(n))/(n (1-p) \lam).$ We use the notation $f(n)
  \stackrel{<}{\sim} g(n)$ to mean that $f(n)$ is asymptotically
  bounded by a function $b h(n)$ where $b$ is a constant and $h(n)
  \sim g(n).$

  Let $n_k = k^a,$ be a subsequence with constant $a$ to be chosen
  later. Let $N_k = \lfloor n_k/p \rfloor.$
  \begin{eqnarray*}
    \prob{\cup_{n=n_k}^{n_{k+1}} \left( c_n \geq y_n \right) } & \leq
    & \prob{c_{n_k} \geq \frac{p (1+\ep)\ln(n_k)}{n_{k+1} (1-p) \lam}}\\
    & \sim & \prob{c_{n_k}^* \geq \frac{p(1+\ep)\ln(n_k)}{n_{k+1}(1-p)
    \lam}} \\
    & = & \prob{\cup_{i=1}^{n_{k-1}} \left( Y_{i,n_k} \geq
    \frac{p(1+\ep)\ln(n_k)}{n_{k+1}(1-p)\lam} \right) } \\
    & \leq & \sum_{i=1}^{n_k - 1} \exp\left( \frac{(N_k - i) p
    (1+\ep)\ln(n_k)}{n_{k+1}(1-p)} \right) \\
    & = & \sum_{j=N_k - n_k +1}^{N_k -1} \left(\frac{1}{n_k}
    \right)^{\frac{j p (1+\ep)}{n_k (1-p) }} \\
    & < & n_k \left( \frac{1}{n_k} \right)^{\frac{(N_k - n_k +1)p
        (1+\ep)}{n_k (1-p)}} \sim \left( \frac{1}{n_k} \right)
        ^{\frac{(N_k/n_k - 1)p(1+\ep)}{(1-p)} -1}\\
    & \stackrel{<}{\sim} & \frac{1}{k^{a \ep}},
  \end{eqnarray*}
where the last line follows by using the fact that $N_k/n_k \rar
p^{-1},$ as $ k \rar \infty.$ Thus, for any $a>1/\epsilon$, we get
  \[
  \sum_{k=0}^{\infty} \prob{\cup_{n=n_k}^{n_{k+1}} \left( c_n \geq
    \frac{p (1+\ep)\ln(n)}{n (1-p) \lam} \right) } < \infty \qquad
    \forall \ep > 0.
  \]
  It follows from the Borel-Cantelli lemma that %
  \[
  \limsup_{n \rar\infty} \frac{ \lam n c_n}{ \ln(n)}
  \leq \frac{p}{1-p} \qquad a.s.
  \]
  To establish the lower bound, we take $y_n = p (1-\ep)\ln(n)/(n (1-p)
  \lam),$ and show that $\prob{c_n^* \leq y_n}$ is summable;
  \[
  \prob{c_n^* \leq y_n} = \prod_{j=N-n+1}^{N-1} (1 - e^{-\lam
    j y_n}).
  \]
  Since $y_n \rar 0,$ and $ny_n \rar \infty$, we have
  \[
  \ln(\prob{c_n^* \leq y_n} \stackrel{<}{\sim} -\frac{(1-p) n}
     {p(1-\ep)\ln(n)},
  \]
  and hence $\prob{c_n^* \leq y_n}$ is summable. This completes the
  proof of the first part.

  Proof for the largest nearest neighbor distance is similar. In the
  proof of the upper bound we take $y_n = p (1+\ep)\ln(n)/(2 n \lam
  (1-p))$ and use the following inequalities.  %
  \begin{eqnarray*}
    \lefteqn{\prob{ d_n^* \geq y_n}}\\
    & = & \prob{\cup_{i=2}^{n-1}((Y_{i-1,n}
      \geq y_n)\cap (Y_{i,n} \geq y_n)) \cup (Y_{1,n} \geq y_n) \cup
      (Y_{n-1,n} \geq y_n)} \\
    & \leq & \sum_{i=2}^{n-1} \prob{Y_{i-1,n} \geq y_n} \prob{Y_{i,n}
      \geq y_n} + \prob{Y_{1,n} \geq y_n} + \prob{Y_{n-1,n} \geq
      y_n}.
  \end{eqnarray*}
  To prove the lower bound set $y_n = p(1-\ep)\ln(n)/(2 n \lam (1-p))$
  and proceed as follows.
  \begin{eqnarray*}
    \lefteqn{ \prob{ d_n^* \leq y_n} }\\
    & =  & \prob{\cap_{i=2}^{n-1}((Y_{i-1,n}
      \leq y_n)\cup (Y_{i,n} \leq y_n)) \cap (Y_{1,n} \leq y_n) \cap
      (Y_{n-1,n} \leq y_n)} \\
    & \leq & \prob{\cap_{i=1}^{\lfloor n/2 \rfloor} ((Y_{2i-1,n} \leq
      y_n) \cup (Y_{2i,n} \leq y_n))} \\
    & \leq & \prod_{i=1}^{\lfloor n/2 \rfloor}(1-e^{-\lam (2N - 4i
      -1)y_n})  .
  \end{eqnarray*}
  Take logarithms and use appropriate Taylor expansions and
  asymptotic equivalences as in the proof of the first part. We
  omit the details. This completes the proof of the Theorem.
  \hfill $\Box$
}

{\bf Acknowledgement: } The authors would like to thank an
anonymous referee for careful reading of the paper and making
many suggestions resulting in significant improvements.

\bibliography{}

\end{document}